\documentclass[12pt, a4paper]{article}
\usepackage{amssymb, amsmath,amsthm }
\usepackage[all]{xy}
\newtheorem{prop}{Proposition}[section]
\newtheorem{defi}[prop]{Definition}

\newtheorem{lem}[prop]{Lemma}
\newtheorem{thm}[prop]{Theorem}

\newtheorem{cor}[prop]{Corollary}

\DeclareMathAlphabet{\mathpzc}{OT1}{pzc}{m}{it}

\DeclareMathOperator{\End}{End}
\DeclareMathOperator{\Hom}{Hom}
\DeclareMathOperator{\Ind}{Ind}
\DeclareMathOperator{\cInd}{c-Ind}

\DeclareMathOperator{\GL}{GL}

\DeclareMathOperator{\Ker}{Ker}

\DeclareMathOperator{\supp}{Supp}

\DeclareMathOperator{\St}{St}
\DeclareMathOperator{\Sp}{Sp}

\newcommand{\cIndu}[3]{\cInd_{#1}^{#2}{#3}}

\newcommand{\Indu}[3]{\Ind_{#1}^{#2}{#3}}

\newcommand{\FF}{\mathcal F}

\newcommand{\Fq}{\mathbf{F}_{q}}
\newcommand{\Fbar}{\overline{\mathbf{F}}_{p}}

\newcommand{\oo}{\mathfrak o}
\newcommand{\oF}{\oo}
\newcommand{\pF}{\mathfrak{p}}

\newcommand{\pif}{\varpi}

\newcommand{\QK}[1]{I_{1}}

\newcommand{\HH}{\mathcal H}

\newcommand{\Qp}{\mathbf Q_p}

\newcommand{\QQ}{\mathbf{Q}}

\newcommand{\Eins}{\boldsymbol 1}

\newcommand{\oK}{\mathcal O}
\title{On the restriction of  representations of $\GL_2(F)$ to a Borel subgroup }
\author{Vytautas Paskunas}
\date{\today.}
\begin{document} 
\maketitle

\tableofcontents
\pagebreak

\section{Introduction}
Let $F$ be a non-Archimedean local field and let $p$ be the residual characteristic of $F$. Let $G=\GL_2(F)$ and let $P$ be a Borel subgroup of $G$. 
In this paper we study the restriction of irreducible $\Fbar$-representations of $G$ to $P$. We show that in a certain sense 
$P$ controls the representation theory of $G$. We then extend our results to smooth $\oK[G]$- modules of finite length  and unitary $K$-Banach 
space representations of $G$, where $\oK$ is the ring of integers of a complete discretely valued field $K$, with residue field $\Fbar$ and uniformizer 
$\varpi_K$. 

The study of smooth irreducible $\Fbar$-representations of $G$ have been initiated by Barthel and Livne in \cite{bl}. They have shown that 
smooth irreducible $\Fbar$-representations of $G$ with central character fall into four classes:
\begin{itemize}
\item[(1)] one-dimensional representations $\chi\circ\det$;
\item[(2)] (irreducible) principal series $\Indu{P}{G}{(\chi_1\otimes\chi_2)}$, with $\chi_1\neq \chi_2$;
\item[(3)] special series $\Sp\otimes\chi\circ \det$;
\item[(4)] supersingular.
\end{itemize}
Here, $\Sp$ is defined by an exact sequence 
$$0\rightarrow \Eins \rightarrow \Indu{P}{G}{\Eins}\rightarrow \Sp\rightarrow 0,$$
and the supersingular representations can be characterised by the fact that they are not subquotients of  $\Indu{P}{G}{\chi}$
for any smooth character $\chi:P\rightarrow \Fbar^{\times}$. Such representations have been classified only in the case when $F=\Qp$, 
by Breuil \cite{breuil}. 
If $F\neq \Qp$ no such classification is known so far although in a joint work with Breuil 
we can show that there are ``a lot  more'' supersingular representations  
than in the case $F=\Qp$.

The main result of this paper can be summed up as follows:
\begin{thm}\label{main} Let $\pi$ and $\pi'$ be 
smooth $\Fbar$-representations of $G$, such that $\pi$ is irreducible with a central character, then the following hold:
\begin{itemize}
\item[(i)] if $\pi$ is in the principal series then $\pi|_P$ is of length $2$, otherwise  $\pi|_P$ is an irreducible representation of $P$;
\item[(ii)] We have 
$$\Hom_P(\Sp, \pi')\cong \Hom_G(\Indu{P}{G}{\Eins}, \pi'),$$
and if $\pi$ is not in the special series then
$$\Hom_{P}(\pi, \pi')\cong \Hom_G(\pi, \pi').$$
\end{itemize}
\end{thm} 
The first part of this Theorem and the second part with $\pi'$ irreducible is due to Berger \cite{berger} 
in the case $F=\Qp$. Berger uses the theory 
of $(\phi, \Gamma)$-modules and the classification of supersingular representations. 
Our proof is completely different and purely representation theoretic. In fact this paper grew out of trying to find a simple representation theoretic reason, 
to explain Berger's results.  Vigneras in \cite{vig} 
has studied the restriction of principal series representation  of split reductive $p$-adic groups to a Borel subgroup. Her results contain the first part 
of the Theorem in the case  $\pi$ is not supersingular and $F$ arbitrary.   

Using the theorem we extend the result to smooth $\oK[G]$ modules of finite length.

\begin{thm} Let $\pi$ and $\pi'$ be smooth $\oK[G]$ modules and suppose that $\pi$ is of finite length, and the irreducible subquotients 
of $\pi$ admit a central character. Let $\phi\in \Hom_{\oK[P]}(\pi, \pi')$ and suppose that $\phi$ is not $G$-equivariant. Let $\tau$ be 
be the maximal submodule of $\pi$, such that $\phi|_{\tau}$ is $G$-equivariant, and let $\sigma$ be an  irreducible $G$-submodule 
of $\pi/\tau$, then
$$\sigma \cong\Sp\otimes \delta\circ \det,$$
 for some smooth character $\delta: F^{\times} \rightarrow \Fbar^{\times}.$ Moreover, choose $v\in \pi$ such that the image $\overline{v}$ in $\sigma$
spans $\sigma^{I_1}$, then $\Pi \phi(v)-\phi(\Pi v)\neq 0$, $\varpi_K(\Pi \phi(v)-\phi(\Pi v))=0$, and 
$$ g(\Pi \phi(v)-\phi(\Pi v))=\delta(\det g) (\Pi \phi(v)-\phi(\Pi v)), \quad \forall g\in G,$$
where $\Pi$ and $I_1$ are defined in \S\ref{not}.
\end{thm}

This criterion implies:
\begin{cor} Let $\Pi_1$ and $\Pi_2$ be unitary $K$-Banach space representations of $G$.  
Let $\lVert\centerdot\lVert_1$ and $\lVert\centerdot \lVert_2$ be $G$-invariant norms 
defining the topology on $\Pi_1$ and $\Pi_2$. 
Set
 $$L_1=\{v\in\Pi_1:\lVert v \lVert_1\le 1\},\quad L_2=\{v\in\Pi_2:\lVert v \lVert_2\le 1\}.$$
 Suppose that $L_1\otimes_{\oK}\Fbar$ is of finite length as 
$\oK[G]$ module and the irreducible subquotients
admit a central character. Moreover, suppose that
 if $\Sp\otimes\delta\circ\det$ is a subquotient of  
$L_1\otimes_{\oK}\Fbar$, then $\delta\circ\det$ is not a subobject of $L_2\otimes_{\oK}\Fbar$, then 
$$\mathcal L_G(\Pi_1, \Pi_2)\cong \mathcal L_P(\Pi_1,\Pi_2),$$
where $\mathcal L(\Pi_1,\Pi_2)$ denotes continuous $K$-linear maps. 
\end{cor}

Moreover, Theorem \ref{main} implies: 

\begin{cor} Let $\Pi$ be a  unitary $K$-Banach space representation of $G$, let $\lVert\centerdot\lVert$ 
be a $G$-invariant norm defining the topology on $\Pi$. Set
$$L=\{v\in \Pi: \lVert v\lVert\le 1\}.$$
Suppose that $L\otimes_{\oK}\Fbar$ is a finite length $\oK[G]$ module, and the 
irreducible subquotients are either supersingular or characters, then 
every closed $P$-invariant subspace of $\Pi$ is also $G$-invariant.
\end{cor}

According to Breuil's $p$-adic Langlands philosophy  a $2$-dimensional $p$-adic representation 
of the absolute Galois group of $F$  
should be related to a unitary $K$-Banach space representation  of $G$, 
see a forthcoming work of Colmez for the case $F=\Qp$, 
where the restriction to a Borel subgroup plays a prominent role. 
However, if $F\neq \Qp$ it is an open problem to 
construct such unitary $K$-Banach space representations of $G$. 
We hope that our results  will help to understand this.

\textit{Acknowledgements.} This paper was written while I was working with
 Chris\-tophe Breuil on a related project. I would like to thank 
Christophe Breuil for his comments and for pointing out some errors in an earlier draft. 
I would like to thank Eike Lau for a stimulating 
discussion, which led to a simplification of proofs in section \S\ref{last}. I would like to thank
Florian Herzig and Marie-France Vigneras, their comments improved the original manuscript.    

\section{Notations}\label{not}
Let $\oF$ be the ring of integers of $F$, $\pF$ the maximal ideal of $\oF$, and let $q$ be the number of elements in the residue field $\oF/\pF$. 
We fix a uniformiser $\pif$ and an embedding $\oF/\pF\hookrightarrow \Fbar$. 
For $\lambda\in \Fq$ we denote the Teichm\"uller lift of $\lambda$ to $\oF$ by $[\lambda]$. Set 
$$\Pi=\begin{pmatrix} 0 & 1\\ \pif & 0\end{pmatrix}, \quad s=\begin{pmatrix} 0 & 1 \\ 1 & 0\end{pmatrix}, \quad t= \begin{pmatrix} \pif & 0 \\ 0 & 1 \end{pmatrix}. $$
Let $P$ be subgroup of upper-triangular matrices in $G$, $T$ the subgroup of diagonal matrices, $K=\GL_2(\oF)$ and 
$$ I=\begin{pmatrix} \oF^{\times} & \oF \\ \pF & \oF^{\times} \end{pmatrix},\quad I_1= \begin{pmatrix} 1+\pF & \oF \\ \pF & 1+\pF \end{pmatrix}, 
\quad K_1=\begin{pmatrix} 1+\pF & \pF \\ \pF & 1+\pF \end{pmatrix}.$$
All the representations in this paper are on $\Fbar$-vector spaces, except for section \S\ref{last}.

\section{Key}\label{key}
In this section we show how to control the action of $s$ on a supersingular representation $\pi$ in terms of the action of $P$. All the hard work 
here is done by Barthel and Livne in \cite{bl}, we just record a consequence of their proof of \cite{bl} Theorem 33.  
 
Let $\sigma$ be an irreducible representation of $K$. Let $\tilde{\sigma}$ be a representation
of $F^{\times}K$ such that $\pif$ acts trivially on $\tilde{\sigma}$ and $\tilde{\sigma}|_K=\sigma$. 
Set $\FF_{\sigma}=\cIndu{F^{\times}K}{G}{\tilde{\sigma}}$ and 
$\HH_{\sigma}=\End_G(\FF_{\sigma})$. It is shown in 
\cite{bl} Proposition 8 that as an algebra $\HH_{\sigma}\cong \Fbar[T]$, for a 
certain $T\in \HH_{\sigma}$, defined in 
\cite{bl} \S3. Fix $\varphi\in \FF_{\sigma}$ such that 
 $\supp \varphi =F^{\times}K$ and $\varphi(1)$ spans $\sigma^{I_1}$. Since $\varphi$ generates $\FF_{\sigma}$ as 
a $G$-representation $T$ is determined by $T\varphi$. 
\begin{lem}\label{T}
 \begin{itemize}
\item[(i)] If $\sigma \cong \psi\circ \det$, for some character $\psi:
\oF^{\times}\rightarrow \Fbar^{\times}$, then
$$ T \varphi = \Pi \varphi +\sum_{\lambda \in \Fq} \begin{pmatrix} 1 & [\lambda]\\ 0 & 1
\end{pmatrix} t \varphi .$$
\item[(ii)] Otherwise, 
$$ T \varphi = \sum_{\lambda \in \Fq} \begin{pmatrix} 1 & [\lambda] \\ 0 & 1
\end{pmatrix} t  \varphi.$$
\end{itemize}
\end{lem}
\begin{proof} In the notation of \cite{bl} this is a calculation of $T([1, e_{\vec{0}}])$. The claim follows from the formula (19) in the proof 
of \cite{bl} Theorem 19. 
 \end{proof}

Let $\pi$ be a supersingular representation of $G$, such that $\pif$ acts trivially. Let $v\in \pi^{I_1}$ and suppose that 
$\langle K\centerdot v\rangle\cong \sigma$. The Frobenius reciprocity gives $\alpha\in \Hom_{G}(\FF_{\sigma}, \pi)$, such that 
$\alpha(\varphi)=v$.
\begin{lem}\label{nilp} There exists an $n\ge 1$ such that $\alpha\circ T^n=0$.
\end{lem}
\begin{proof} Now $\Hom_G(\FF_{\sigma}, \pi)$ is naturally a right $\HH_{\sigma}$-module; let $M=\langle \alpha\centerdot \HH_{\sigma}\rangle$ be 
an $\HH_{\sigma}$-submodule of $\Hom_G(\FF_{\sigma}, \pi)$ generated by $\alpha$. The proof of \cite{bl} Proposition 32 implies that $\dim_{\Fbar} M$ 
is finite. Let $\overline{T}$ be the image of $T$ in $\End_{\Fbar}(M)$ and let $m(X)$ be the minimal polynomial of $\overline{T}$. Let   
$\lambda\in \Fbar$ be such that $m(\lambda)=0$, then we may write
$m(X)= (X-\lambda)h(X)$. Since $m(X)$ is minimal the composition 
$$ h(T)(\FF_{\sigma}) \rightarrow \FF_{\sigma}\rightarrow \pi$$    
is non-zero. According to \cite{bl} Theorem 19, $\FF_{\sigma}$ is a free $\HH_{\sigma}$ module, hence $h(T)$ is an injection and so 
$h(T)(\FF_{\sigma})$ is isomorphic to $\FF_{\sigma}$. This implies that $\pi$ is a quotient of $\FF_{\sigma}/(T-\lambda)$. Since $\pi$ is supersingular
\cite{bl} Corollary 36 implies that $\lambda=0$, and hence $m(X)=X^n$, for some $n\ge 1$. 
\end{proof}

\begin{cor}\label{ex} Let $\pi$ be a supersingular representation, such that $\pif$ acts trivially. Let $v\in \pi^{I_1}$ be such that
 $\langle K\centerdot v\rangle$ is an irreducible representation of $K$. Set $v_0=v$ and for $i\ge 0$ set
$$ v_{i+1}= \sum_{\lambda\in \Fq} \begin{pmatrix} 1 & [\lambda]\\ 0 & 1\end{pmatrix} t v_i.$$ 
Then $v_i\in \pi^{I_1}$ for all $i\ge 1$ and there exists an $n\ge 1$, such that $v_n=0$.
\end{cor}
\begin{proof} Set $\sigma= \langle K\centerdot v\rangle$. If $\sigma$ is not a character then Lemma \ref{T} (ii) implies that 
$v_i= (\alpha \circ T^i)(\varphi)$, for all $i\ge 0$ in particular $I_1$ acts trivially on $v_i$ and 
the statement follows from Lemma \ref{nilp}. If $\sigma$ is a character, then after 
twisting we may assume that $\sigma=\Eins$.  Since $I$ acts trivially on $\Pi v_0$  
the space $\langle K\centerdot (\Pi v_0)\rangle$ is a quotient of $\Indu{I}{K}{\Eins}$. Now 
$$v_1= \sum_{\lambda\in \Fq} \begin{pmatrix} 1 & [\lambda] \\ 0 & 1\end{pmatrix} s ( \Pi v_0).$$ 
If $v_1=0$ then we are done. If  $v_1\neq 0$ then \cite{p} (3.1.7) and (3.1.8) imply that $\langle K\centerdot v_1\rangle\cong \St$, where 
$\St$ is the inflation of the Steinberg representation of $\GL_2(\Fq)$. We may apply the previous part to $v_1$.
\end{proof}

\begin{lem}\label{s} Let $\pi$ be a smooth representation of $G$ and let $v\in \pi^{I_1}$. Suppose that
$$ \sum_{\lambda\in \Fq} \begin{pmatrix} 1 & [\lambda]\\ 0 & 1\end{pmatrix} t v = 0.$$
Then 
$$ s v= - \sum_{\lambda\in \Fq^{\times}} \begin{pmatrix} - \pif [\lambda^{-1}] & 1 \\ 0 & \pif^{-1}[\lambda]\end{pmatrix} v.$$
\end{lem}
\begin{proof} We may rewrite 
$$ v= -\sum_{\lambda\in \Fq^{\times}} t^{-1} \begin{pmatrix} 1 & [\lambda] \\ 0 & 1\end{pmatrix} t v =- \sum_{\lambda\in \Fq^{\times}}
 \begin{pmatrix}  1 & \pif^{-1} [\lambda] \\ 0 & 1 \end{pmatrix} v.$$
If $\beta\in F^{\times}$ then 
\begin{equation}\label{trix}
\begin{pmatrix} 0 & 1 \\ 1 & 0\end{pmatrix} \begin{pmatrix} 1 & \beta \\ 0 & 1\end{pmatrix}= 
\begin{pmatrix} -\beta^{-1} & 1 \\ 0 & \beta \end{pmatrix}\begin{pmatrix} 1 & 0 \\ \beta^{-1} & 1\end{pmatrix}.
\end{equation}
Since $v$ is invariant by $\begin{pmatrix} 1 & 0 \\ \pF & 1 \end{pmatrix}$ the matrix identity above implies the result.
\end{proof}

Since $G= P I_1 \cup P s I_1$, we will use the Lemma above to show that the action of $P$ on $\pi$ already ``contains all the information''
about the action of $G$ on $\pi$.

\section{Supersingular representations}
In this section we study the restriction  of a supersingular representations of $G$ to a Borel subgroup.

\begin{lem}\label{next} Let $\pi$ be a smooth representation of $G$ and let  $v\in \pi^{I_1}$ be non-zero,
and such that $I$ acts on $v$ via  a character $\chi$,  
then there exists $j$, such that  $0\le j \le q-1$ and if we let 
$$w=\sum_{\lambda\in\Fq} \lambda^j \begin{pmatrix} 1 & [\lambda] \\ 0 & 1 \end{pmatrix} t v,$$
then $w\in \pi^{I_1}$ and $\langle K\centerdot w\rangle$ is an irreducible representation of $K$. 
\end{lem}

\begin{proof} Set $\tau=\langle K \centerdot (\Pi v)\rangle$.  For $0\le j \le q-1$ set 
$$w_j=\sum_{\lambda\in\Fq} \lambda^j \begin{pmatrix} 1 & [\lambda] \\ 0 & 1 \end{pmatrix} s(\Pi v)=
\sum_{\lambda\in\Fq} \lambda^j \begin{pmatrix} 1 & [\lambda] \\ 0 & 1 \end{pmatrix} t v.$$
The set $\{\Pi v, w_j: 0\le j \le q-1\}$ spans $\tau$. 

If $w_0=0$ then Lemma \ref{s} implies that 
$$ \Pi v= -\sum_{\lambda\in \Fq^{\times}}\begin{pmatrix}-\pif[\lambda^{-1}] & 1 \\ 0 & [\lambda]\end{pmatrix} v=  
-\sum_{\mu\in \Fq^{\times}}\chi(\begin{pmatrix} -[\mu] & 0 \\ 0 & [\mu^{-1}]\end{pmatrix})\begin{pmatrix}1 & [\mu] \\ 0 & 1\end{pmatrix} t v.$$
Since 
$$ \chi(\begin{pmatrix} [\mu] & 0 \\ 0 & [\mu^{-1}]\end{pmatrix})=\mu^r, \quad \forall \mu\in \Fq^{\times}$$
for some $0\le r< q-1$, we obtain that $\tau$ is spanned by the set $\{w_j: 1\le j\le q-1\}$. 
Let $\sigma$ be a $K$-irreducible 
subrepresentation of $\tau$. 
The space $\sigma^{I_1}$ is one dimensional, so $I$ acts on $\sigma^{I_1}$ by a character. However, one 
may verify that the group
$$\biggl \{ \begin{pmatrix} [\lambda] & 0 \\ 0 & 1\end{pmatrix}: \lambda\in \Fq^{\times} \biggr \}$$
acts on the set $w_j$ for $1\le j \le q-1$ by distinct characters, hence $\sigma^{I_1}$ is spanned by $w_j$ for some $1\le j\le q-1$.

Suppose that $w_0\neq 0$. If $w_0$ and $\Pi v$ are linearly independent
 then the natural map 
$\Indu{I}{K}{\chi^s}\rightarrow \tau$ is an injection, since it induces an injection 
on $(\Indu{I}{K}{\chi^s})^{I_1}$. It follows from \cite{p} (3.1.5) that 
$\langle K\centerdot w_0\rangle$ is an irreducible representation of $K$. 
If $w_0$ and $\Pi v$ are not linearly independent then $\chi=\chi^s$. It follows from \cite{p} (3.1.8) 
that $\langle K\centerdot w_0\rangle$ is a isomorphic to a twist of the Steinberg representation by a 
character.
 
\end{proof}  

\begin{prop}\label{give} Let $\pi$ be a smooth representation of $G$, and let $w$ be a non-zero vector in $\pi$. Then there exists 
a non-zero $v\in \langle P\centerdot w\rangle \cap \pi^{I_1}$ such that $\langle K\centerdot v \rangle$ is an irreducible 
representation of $K$. 
\end{prop}
\begin{proof} Since $\pi$ is smooth there exists $k\ge 0$ such that $w$ is fixed by 
$\begin{pmatrix} 1 & 0 \\ \pF^{k+1} & 1\end{pmatrix}$. Then $w_1= t^k w$
is fixed by $\begin{pmatrix} 1 & 0 \\ \pF & 1\end{pmatrix}$. Iwahori decomposition gives us 
$$ I_1= \begin{pmatrix} 1+\pF & \oF \\ 0 & 1+\pF\end{pmatrix}  \begin{pmatrix} 1 & 0 \\ \pF & 1\end{pmatrix}.$$
Hence, $\tau:=\langle I_1\centerdot w_1\rangle= \langle (I_1\cap P)\centerdot w_1\rangle$. Since $I_1$ is 
a pro-$p$ group, we have $\tau^{I_1}\neq 0$, and hence $\langle P\centerdot w\rangle \cap \pi^{I_1} \neq 0$. 
Let $w_2\in \langle P\centerdot w\rangle \cap \pi^{I_1} \neq 0$ be non-zero. Since $|I/ I_1|$ is prime to $p$, 
there exists a smooth character $\chi:I \rightarrow \Fbar^{\times}$ such that
 $$w_3=\sum_{\lambda, \mu\in \Fq^{\times}} \chi(\begin{pmatrix} [\lambda^{-1}] & 0 \\ 0 & [\mu^{-1}] \end{pmatrix})  
\begin{pmatrix} [\lambda] & 0 \\ 0 & [\mu] \end{pmatrix} w_2$$
is non-zero. Lemma \ref{next} applied to $w_3$ gives the required vector.  

\end{proof}      

\begin{thm}\label{irred} Let $\pi$ be  supersingular, then $\pi|_P$ is an irreducible representation of $P$. 
\end{thm}
\begin{proof} Let $w\in \pi$ be non-zero. According to Proposition \ref{give}
there exists  a non-zero 
$v\in \langle P\centerdot w\rangle \cap \pi^{I_1}$, such that $\sigma:=\langle K\centerdot v\rangle $ is an irreducible 
representation of $K$. Corollary \ref{ex} implies that there exists a non-zero $v'\in \pi^{I_1} \cap \langle P\centerdot v\rangle$ such that 
$$\sum_{\lambda\in \Fq} \begin{pmatrix} 1 & [\lambda] \\ 0 & 1\end{pmatrix} t v'=0.$$
According to Lemma \ref{s} $s v'\in \langle P\centerdot v'\rangle$. Since $G=P I_1 \cup P s I_1$ and $\pi$ is irreducible $G$-representation we have 
$$ \pi= \langle G\centerdot v'\rangle= \langle P\centerdot v'\rangle \subseteq \langle P\centerdot w\rangle.$$
Hence, $\pi=\langle P\centerdot w\rangle$ for all $w\in \pi$, and so $\pi|_P$ is irreducible.   
\end{proof} 

\begin{thm}\label{restP} Let $\pi$ and $\pi'$ be smooth  representations of $G$, such that $\pi$ is supersingular, then
$$\Hom_P(\pi, \pi')\cong \Hom_G(\pi,\pi').$$
\end{thm} 

\begin{proof} Let $\phi\in \Hom_P(\pi,\pi')$ be non-zero. Choose $v\in \pi^{I_1}$ such that 
$\langle K\centerdot v\rangle$ is an irreducible representation of $K$. 
Since by Theorem \ref{irred} $\pi|_P$ is  irreducible $\phi$ is an injection and hence $\phi(v)\neq 0$.
Since $v$ is fixed by $I_1$ and $\phi$ is $P$-equivariant, we have that $\phi(v)$ is fixed by $I_1\cap P$.
 Since $\pi'$ is smooth there exists an integer $k\ge 1$  such that  
$\phi(v)$ is fixed by $\begin{pmatrix} 1 & 0 \\ \pF^k & 1 \end{pmatrix}$. Suppose that $k> 1$. 
Lemma \ref{next} implies that there exists $j$, such that $0\le j\le q-1$ and if we set
$$v_1=\sum_{\lambda\in\Fq} \lambda^j \begin{pmatrix} 1 & [\lambda] \\ 0 & 1 \end{pmatrix} t v,$$
then $v_1\in \pi^{I_1}$ and $\langle K\centerdot v_1\rangle$ is an irreducible representation of $K$. Since $\phi$ is $P$-equivariant, $\phi(v_1)$
is fixed by $I_1\cap P$ and 
$$\phi(v_1)=\sum_{\lambda\in\Fq} \lambda^j \begin{pmatrix} 1 & [\lambda] \\ 0 & 1 \end{pmatrix} t \phi(v).$$
If $\alpha\in \oF$ and $\beta\in \pF$ then 
$$\begin{pmatrix} 1 & 0 \\ \beta & 1\end{pmatrix} \begin{pmatrix} 1 & \alpha \\ 0 & 1\end{pmatrix}=
 \begin{pmatrix} 1 & \alpha(1+\alpha\beta)^{-1} \\ 0 & 1\end{pmatrix} \begin{pmatrix} (1+\alpha\beta)^{-1} & 0 \\ 
\beta & 1+\alpha\beta\end{pmatrix}.$$
This matrix identity coupled with 
$$t^{-1}\begin{pmatrix} 1 & 0 \\ \pF^{k-1} & 1 \end{pmatrix} t= \begin{pmatrix} 1 & 0 \\ \pF^k & 1\end{pmatrix},$$
implies that $\phi(v_1)$ is fixed by $\begin{pmatrix} 1 & 0 \\ \pF^{k-1} & 1 \end{pmatrix}$. 
By repeating the argument we obtain $w\in \pi^{I_1}$ such that $\langle K\centerdot w\rangle$ is an irreducible representation of $K$ 
and $\phi(w)$ is fixed by 
$\begin{pmatrix} 1 & 0 \\ \pF & 1\end{pmatrix}$. Iwahori decomposition implies that $\phi(w)$ is fixed by $I_1$. Set $v_0=w$ and for $i\ge 0$,
$$ v_{i+1}= \sum_{\lambda\in \Fq} \begin{pmatrix} 1 & [\lambda]\\ 0 & 1\end{pmatrix} t v_i.$$ 
Since $v_i$ are fixed by $I_1$, $\phi(v_i)$ are fixed by $I_1\cap P$. Moreover, 
$$ \phi(v_{i+1})= \sum_{\lambda\in \Fq} \begin{pmatrix} 1 & [\lambda]\\ 0 & 1\end{pmatrix} t \phi(v_i).$$  
Since $\phi(v_0)$ is fixed by $I_1$, the  argument used above implies that $\phi(v_{i+1})$ are fixed by 
$\begin{pmatrix} 1 & 0 \\ \pF & 1\end{pmatrix}$ and hence  fixed by $I_1$. 
 Corollary \ref{ex} implies that $v_n=0$ for some $n\ge 1$.
Let $m$ be the smallest integer such that $v_m=0$ and set $v'=v_{m-1}$. Then  $v'\in \pi^{I_1}$,  $\phi(v')\in (\pi')^{I_1}$ and 
$$ \sum_{\lambda\in \Fq}\begin{pmatrix} 1 & [\lambda]\\ 0 & 1\end{pmatrix} t v'=0, \quad 
 \sum_{\lambda\in \Fq}\begin{pmatrix} 1 & [\lambda]\\ 0 & 1\end{pmatrix} t \phi(v')=0.$$ 
Lemma \ref{s} applied to $v'$ and $\phi(v')$ implies that 
\begin{equation}\notag
\begin{split}
 \phi(s v')=&- \phi \biggl(\sum_{\lambda\in \Fq^{\times}} \begin{pmatrix} - \pif [\lambda^{-1}] & 1 \\ 0 & \pif^{-1}[\lambda]\end{pmatrix} v'\biggr )\\
           =&- \sum_{\lambda\in \Fq^{\times}} \begin{pmatrix} - \pif [\lambda^{-1}] & 1 \\ 0 & \pif^{-1}[\lambda]\end{pmatrix} \phi(v')=s \phi(v')
\end{split}
\end{equation}

Since $G=P I_1\cup P s I_1$ this implies that $\phi( \pi(g) v')= \pi'(g)\phi(v')$, for all $g\in G$. Since $\pi$ is irreducible 
$\pi=\langle G\centerdot v'\rangle$ and this implies that $\phi$ is $G$-equivariant.
\end{proof} 

\section{Non-supersingular representations}
Let $\chi: T\rightarrow \Fbar^{\times}$ be a smooth character. We consider it as a character of $P$, via 
$P\rightarrow P/U\cong T$. We define a smooth representation $\kappa_{\chi}$ of $P$ by the short exact sequence:
\begin{equation}\label{kappa}
0 \rightarrow \kappa_{\chi} \rightarrow \Indu{P}{G}{\chi}\rightarrow \chi\rightarrow 0
\end{equation}
where the map on the right is given by the evaluation at the identity. 
The representation $\kappa_{\chi}$ is absolutely irreducible by \cite{vig} Theoreme 5. 
If $\chi=\psi\circ\det$ for some smooth character $\psi:F^{\times}\rightarrow \Fbar^{\times}$ then the sequence splits as a $P$-representation and we obtain
$$\Sp\otimes \psi\circ\det|_P \cong \kappa_{\chi}.$$
 
\begin{lem}\label{char} Let  $\pi$ be a smooth representation of $G$. Suppose that 
$\Hom_P(\chi, \pi)\neq 0$
 then $\chi$ extends uniquely to a character of $G$, and 
$$\Hom_P(\chi, \pi)\cong \Hom_G(\chi, \pi).$$
\end{lem}
\begin{proof} Let $\phi\in \Hom_P(\chi, \pi)$ be non-zero, and let $v$ be a basis vector of the underlying vector space of $\chi$. Since $\pi$ is smooth 
there exists $k\ge 1$ such that $\phi(v)$ is fixed by $\begin{pmatrix} 1 & 0\\ \pF^k & 1 \end{pmatrix}$. Since 
$t \phi(v) = \chi(t) \phi(v)$, we obtain that $\phi(v)$ is  fixed by $\begin{pmatrix} 1 & 0 \\ \pF^{k-1} & 1\end{pmatrix}$, and by repeating this we obtain that 
$\phi(v)$ is fixed by $sUs$. Now $sUs$ and $P$ generate $ G$. This implies the claim.
\end{proof}    

\begin{cor}\label{princinj} 
 Let $\pi'$ be a smooth representation of $G$. Suppose that $\chi\neq \chi^s$ and let $\phi\in \Hom_P(\Indu{P}{G}{\chi}, \pi')$ be non-zero, then 
$\phi$ is an injection.  
\end{cor}

\begin{proof} Lemma 
\ref{char} implies that $\Hom_P(\chi, \Indu{P}{G}{\chi})=0$. Hence the sequence \eqref{kappa} 
cannot split.  
 So if $\Ker \phi\neq 0$ then $\Ker \phi$ contains $\kappa_{\chi}$. Hence, $\phi$ induces a homomorphism $\bar{\phi}\in \Hom_P(\chi, \pi')$. 
Lemma \ref{char} implies that $\bar{\phi}=0$ and hence $\phi=0$.
\end{proof} 

\begin{cor}\label{pichi} Suppose that $\chi\neq \chi^s$ then 
$$\Hom_P(\Indu{P}{G}{\chi}, \Indu{P}{G}{\chi})\cong \Hom_G(\Indu{P}{G}{\chi}, \Indu{P}{G}{\chi}).$$
\end{cor}
\begin{proof} Suppose that $\phi_1, \phi_2\in \Hom_P(\Indu{P}{G}{\chi}, \Indu{P}{G}{\chi})$
are non-zero, then by Corollary \ref{princinj} the restriction  of $\phi_1$ and $\phi_2$ to $\kappa_{\chi}$ induces non-zero homomorphisms in 
$\Hom_P(\kappa_{\chi}, \kappa_{\chi})$. Since $\kappa_{\chi}$ is absolutely irreducible 
this implies that there exists a scalar $\lambda\in \Fbar^{\times}$ such that 
the restriction of $\phi_1-\lambda \phi_2$ to $\kappa_{\chi}$ is zero.  
Now $\phi_1-\lambda\phi_2 \in \Hom_P(\Indu{P}{G}{\chi}, \Indu{P}{G}{\chi})$ and is not an injection,
hence by Corollary \ref{princinj} it must be equal to zero.
\end{proof}
  
\begin{thm}\label{mapkappa} Let $\pi$ be a smooth representation of $G$, then the restriction to $\kappa_{\chi}$ induces an isomorphism
$$ \iota:\Hom_G(\Indu{P}{G}{\chi}, \pi)\cong\Hom_P(\kappa_{\chi}, \pi).$$
\end{thm}
\begin{proof} If $\chi\neq \chi^s$ then the injectivity of 
$\iota$ is given by Corollary \ref{princinj}. If $\chi= \chi^s$ then the injectivity follows 
from Lemma \ref{char} and  \cite{bl} Theorem 30(1)(b). We are going to show that $\iota$ is surjective.

Let $\varphi_1\in \Indu{P}{G}{\chi}$ be an 
$I_1$ invariant function such that $\supp \varphi_1= P I_1$ and $\varphi_1(1)=1$. Set 
$$ \varphi_2=\sum_{\lambda\in \Fq} \begin{pmatrix} 1 & [\lambda]\\ 0 & 1\end{pmatrix} s \varphi_1.$$
Then $\{\varphi_1, \varphi_2\}$ is a basis of $(\Indu{P}{G}{\chi})^{I_1}$. Since $G=PK$ we have 
$$(\Indu{P}{G}{\chi})^{K_1}\cong \Indu{I}{K}{\chi},$$
as a representation of $K$, and hence $\sigma= \langle K\centerdot \varphi_2\rangle$ is an irreducible representation of $K$, which is not a character.
We let $F^{\times}$ act on $\sigma$ via $\chi$. Frobenius reciprocity gives us a map 
$$ \alpha: \cIndu{F^{\times}K}{G}{\sigma} \rightarrow \Indu{P}{G}{\chi}.$$
It follow from \cite{bl} Theorem 30 (3) that there exists $\lambda\in \Fbar^{\times}$, determined by $\chi$, such that $\alpha$ induces an isomorphism 
$$ \cIndu{F^{\times}K}{G}{\sigma}/(T-\lambda) \cong  \Indu{P}{G}{\chi},$$
where $T\in \End_G(\cIndu{F^{\times}K}{G}{\sigma})$ is as in section \S\ref{key}. Lemma 
\ref{T} implies that 
$$\varphi_2=\lambda^{-1}  (\sum_{\mu \in \Fq} \begin{pmatrix} 1 & [\mu]\\ 0 & 1\end{pmatrix} t \varphi_2).$$

Let $\psi\in \Hom_P(\kappa_\chi, \pi)$ be non-zero. Since $\supp \varphi_2= P s I_1$ we have $\varphi_2(1)=0$ and hence $\varphi_2\in \kappa_{\chi}$. 
Since $\kappa_{\chi}$ is irreducible $\psi(\varphi_2)\neq 0$ and the $P$-equivariance of $\psi$ gives:
\begin{equation}\label{rel}
\psi(\varphi_2)= \lambda^{-1} (\sum_{\mu \in \Fq} \begin{pmatrix} 1 & [\mu]\\ 0 & 1\end{pmatrix} t \psi(\varphi_2)).  
\end{equation}
This equality coupled with the  argument used in the proof of \ref{restP} implies that $\psi(\varphi_2)$ is fixed by 
$\begin{pmatrix} 1 & 0 \\ \pF & 1\end{pmatrix}$. Since $\psi$ is $P$-equivariant $\psi(\varphi_2)$ is fixed by $I_1\cap P$. 
The Iwahori decomposition implies that $\psi(\varphi_2)$ is fixed by  $I_1$. 

So $I_1$ fixes $\Pi \psi(\varphi_2)$ and $I$ acts on $\Pi \psi(\varphi_2)$ via the character $\chi$. Hence 
$\langle K\centerdot \Pi \psi(\varphi_2)\rangle$ is a quotient of $\Indu{I}{K}{\chi}$. Now 

$$ \sum_{\mu\in \Fq} \begin{pmatrix} 1 & [\mu] \\ 0 & 1 \end{pmatrix} s (\Pi \psi(\varphi_2))= 
\psi\biggl(\sum_{\mu\in \Fq} \begin{pmatrix} 
1 & [\mu] \\ 0 & 1 \end{pmatrix} t\varphi_2 \biggr)=\lambda \psi(\varphi_2)\neq 0.$$  

If $\chi|_{T\cap K} \neq \chi^s|_{T\cap K}$ then this 
implies that $\langle K\centerdot \Pi \psi(\varphi_2) \rangle \cong \Indu{I}{K}{\chi}$, and hence 
$\langle K\centerdot \psi(\varphi_2)\rangle \cong \sigma$. 
If $\chi|_{T\cap K}= \psi\circ \det$ for some $\psi:\oF^{\times}\rightarrow \Fbar^{\times}$ then the 
above equality implies that if $\Pi \psi(\varphi_2)$  and $\psi(\varphi_2)$ are linearly independent then 
$$\langle K\centerdot \Pi \psi(\varphi_2)\rangle \cong \Indu{I}{K}{\chi},$$
otherwise   
$$\langle K\centerdot \Pi \psi(\varphi_2)\rangle\cong \St\otimes \psi\circ \det,$$ 
where $\St$ is the lift to $K$ of Steinberg representation of 
$\GL_2(\Fq)$. In both cases we obtain that 
$\langle K\centerdot \psi(\varphi_2)\rangle \cong \St\otimes \psi\circ \det\cong\sigma$.
Hence, $\langle G\centerdot \psi(\varphi_2)\rangle$ is a quotient of $\cIndu{F^{\times}K}{G}{\sigma}$. 
The equation \eqref{rel} and Lemma \ref{T} implies that 
$\langle G\centerdot \psi(\varphi_2)\rangle$ is a quotient of 
$$\cIndu{F^{\times}K}{G}{\sigma}/ (T-\lambda)\cong \Indu{P}{G}{\chi}.$$
Hence, $\iota$ is also surjective. 
\end{proof}

\begin{cor}\label{prince} Suppose that $\chi\neq \chi^s$ and let $\pi$ be a smooth representation of $G$ then 
$$\Hom_G(\Indu{P}{G}{\chi}, \pi)\cong \Hom_P(\Indu{P}{G}{\chi}, \pi).$$
\end{cor}
\begin{proof} Let $\psi\in \Hom_P(\Indu{P}{G}{\chi}, \pi)$ be non-zero. It follows 
from Corollary \ref{princinj} that the composition 
$$ \Indu{P}{G}{\chi}\rightarrow \pi\rightarrow \pi/\langle G\centerdot \psi(\kappa_{\chi})\rangle$$
is zero. Hence the image of $\psi$ is contained in $\langle G\centerdot \psi(\kappa_{\chi})\rangle$. It follows from Theorem \ref{mapkappa}
applied to $\pi=\langle G\centerdot \psi(\kappa_{\chi})\rangle$ and the irreducibility of $\Indu{P}{G}{\chi}$ that $\Indu{P}{G}{\chi}$
is isomorphic to $\langle G\centerdot \psi(\kappa_{\chi})\rangle$ as a $G$-representation. The $G$-equivariance of $\psi$ follows from Corollary \ref{pichi}.
\end{proof}

\begin{cor}\label{sp} Let $\pi$ be a smooth representation of $G$, then 
$$ \Hom_P(\Sp, \pi)\cong \Hom_G(\Indu{P}{G}{\Eins}, \pi).$$
\end{cor}
Note that, $\Hom_G(\Sp, \Indu{P}{G}{\Eins})=0$, but $\Hom_G(\Indu{P}{G}{\Eins}, \Indu{P}{G}{\Eins})\neq 0$, so the above result cannot be improved.

\section{Applications}\label{last}
Let $K$ be a complete discrete valuation field, $\oK$ the ring of 
integers and $\varpi_K$ a uniformizer, and we assume that $\oK/\varpi_K \oK\cong \Fbar$.
 We will extend the results of previous sections to smooth $\oK[G]$ modules of finite length, and, 
after passing to the limit, to unitary $K$-Banach space representations of $G$. 

\begin{thm} Let $\pi$ and $\pi'$ be smooth $\oK[G]$ modules and suppose that $\pi$ is of finite length, and the irreducible subquotients 
of $\pi$ admit a central character. Let $\phi\in \Hom_{\oK[P]}(\pi, \pi')$ and suppose that $\phi$ is not $G$-equivariant. Let $\tau$ be 
be the maximal submodule of $\pi$, such that $\phi|_{\tau}$ is $G$-equivariant, and let $\sigma$ be an  irreducible $G$-submodule 
of $\pi/\tau$, then
$$\sigma \cong\Sp\otimes \delta\circ \det,$$
 for some smooth character $\delta: F^{\times} \rightarrow \Fbar^{\times}.$ Moreover, choose $v\in \pi$ such that the image $\overline{v}$ in $\sigma$
spans $\sigma^{I_1}$, then $\Pi \phi(v)-\phi(\Pi v)\neq 0$, $\varpi_K(\Pi \phi(v)-\phi(\Pi v))=0$, and 
$$ g(\Pi \phi(v)-\phi(\Pi v))=\delta(\det g) (\Pi \phi(v)-\phi(\Pi v)), \quad \forall g\in G.$$
 \end{thm}
   
\begin{proof} We denote by $\Indu{1}{G}{\pi'}$ the space of smooth 
functions from $G$ to the underlying $\oK$ module of $\pi'$, equipped with the $G$ action 
via right translations. Let $\alpha: \pi \rightarrow \Indu{1}{G}{\pi'}$ be a $P$-equivariant map, given by
$$ [\alpha(w)](g)= g\phi(w)-\phi(g w), \quad \forall w\in \pi, \forall g\in G.$$
Then $\tau=\Ker \alpha$. Hence $\alpha$ induces a $P$-equivariant map
$$ \overline{\alpha}: \sigma\rightarrow \Indu{1}{G}{\pi'}.$$
Suppose that $\overline{\alpha}$ is $G$-equivariant, then
$$[ g^{-1} \alpha(g v)](1)=[ g^{-1} \overline{\alpha}(g \overline{v})](1)= [\overline{\alpha}(\overline{v})](1)=[\alpha(v)](1)=0.$$
Hence, $g \phi(v)=\phi(g v)$, for all $g\in G$. So the maximality of $\tau$ implies that 
$\overline{\alpha}$ is not $G$-equivariant.
Hence Theorem \ref{restP}, Lemma \ref{char}, Corollaries \ref{prince} and  \ref{sp} 
imply that 
$$ \sigma\cong \Sp\otimes \delta\circ \det$$
for some smooth character $\delta: F^{\times}\rightarrow \Fbar^{\times}$, and 
$$ \langle G\centerdot \alpha(v) \rangle \cong \Indu{P}{G}{\Eins}\otimes\delta\circ\det.$$  
After twisting we may assume that $\delta$ is the trivial character. 
It follows from \cite{bl} Theorem 30(1)(b) that
$$\Hom_G(\Indu{P}{G}{\Eins}, \Indu{P}{G}{\Eins})\cong \Fbar.$$
Corollary \ref{sp} applied to $\pi=\Indu{P}{G}{\Eins}$ implies that $\overline{\alpha}(\overline{v})$
is a scalar multiple of the function denoted $\varphi_2$ in the proof of Theorem 
\ref{mapkappa}. By construction $\alpha(v)=\overline{\alpha}(\overline{v})$. Hence, 
$\alpha(v)$ is fixed by $I_1$ and  
 $\Pi \alpha(v) +\alpha(v)$ spans the trivial subrepresentation of $G$. In particular,
$$ [\Pi\alpha(v)](1)+[\alpha(v)](1)= [h\Pi\alpha(v)](1)+[h\alpha(v)](1), \quad \forall h\in P.$$
Since $\phi$ is $P$-equivariant we obtain:
$$\Pi\phi(v)-\phi(\Pi v)= h(\Pi \phi(v)-\phi(\Pi v)), \quad \forall h\in P.$$
Suppose that $\Pi\phi(v)=\phi(\Pi v)$. Since $\alpha(v)$ is $I_1$-invariant we obtain
$$ h\Pi u  \phi(v)-\phi(h\Pi u v )= [u \alpha(v)](h\Pi)= [\alpha(v)](h\Pi)= h(\Pi\phi(v)-\phi(\Pi v))=0,$$
for all $h\in P$ and $ u\in I_1.$ And
 $$ h u  \phi(v)-\phi(h u v )= [u \alpha(v)](h)= [\alpha(v)](h)=0, \quad \forall u\in I_1, \forall h\in P.$$
Since $G=PI_1\cup P \Pi I_1$, we obtain that $g\phi(v)=\phi(gv)$, for all $g\in G$, but this  
contradicts the 
maximality of $\tau$. 
So $\Pi \phi(v)-\phi(\Pi v)\neq 0$. Since $\sigma$ is irreducible $\varpi_K \overline{v}=0$, and hence
$$ [\varpi_K \alpha(v)](\Pi)=\varpi_K(  \Pi \phi(v)-\phi(\Pi v))=0,$$   
so $\oK(  \Pi \phi(v)-\phi(\Pi v))=\Fbar(\Pi \phi(v)-\phi(\Pi v))$. Lemma 
\ref{char} implies that $G$ acts trivially on  $\Pi \phi(v)-\phi(\Pi v)$.
\end{proof}

\begin{cor}\label{im} Let $\pi$ and by $\pi'$ be as above and suppose that if $\Sp\otimes\delta\circ \det$ is a subquotient 
of $\pi$ then $\delta\circ \det$ is not a subobject of $\pi'$ then 
$$\Hom_G(\pi,\pi')\cong \Hom_P(\pi, \pi').$$
\end{cor}

\begin{defi} A unitary $K$-Banach space representation $\Pi$ of $G$ is a $K$-Banach space $\Pi$ 
equipped with a $K$-linear action of $G$, such that the map $G\times \Pi\rightarrow \Pi$, 
$(g,v)\mapsto gv$  is continuous and such that the topology
on $\Pi$ is given by a $G$-invariant norm. 
\end{defi}

\begin{cor} Let $\Pi_1$ and $\Pi_2$ be unitary $K$-Banach space representations of $G$.  
Let $\lVert\centerdot\lVert_1$ and $\lVert\centerdot \lVert_2$ be $G$-invariant norms 
defining the topology on $\Pi_1$ and $\Pi_2$. 
Set
 $$L_1=\{v\in\Pi_1:\lVert v \lVert_1\le 1\},\quad L_2=\{v\in\Pi_2:\lVert v \lVert_2\le 1\}.$$
 Suppose that $L_1\otimes_{\oK}\Fbar$ is of finite length as 
$\oK[G]$ module and the irreducible subquotients
admit a central character. Moreover, suppose that
 if $\Sp\otimes\delta\circ\det$ is a subquotient of  
$L_1\otimes_{\oK}\Fbar$, then $\delta\circ\det$ is not a subobject of $L_2\otimes_{\oK}\Fbar$, then 
$$\mathcal L_G(\Pi_1, \Pi_2)\cong \mathcal L_P(\Pi_1,\Pi_2),$$
where $\mathcal L(\Pi_1,\Pi_2)$ denotes continuous $K$-linear maps. 
\end{cor}

\begin{proof} Corollary \ref{im} implies that for all $k\ge 1$ we have 
$$ \Hom_G(L_1/\varpi_K^k L_1, L_2/\varpi_K^{k} L_2)\cong  \Hom_P(L_1/\varpi_K^k L_1, L_2/\varpi_K^{k} L_2).$$
Since $\Hom_{\oK}(L_1/\varpi_K^k L_1, L_2/\varpi_K^{k} L_2)\cong \Hom_{\oK}(L_1,L_2/\varpi_K^{k} L_2)$ 
by passing to the limit we obtain:
$$ \Hom_G( L_1, L_2)\cong  \Hom_P(L_1, L_2).$$
It follows from \cite{schn} Proposition 3.1 that 
$$\mathcal L(\Pi_1, \Pi_2)\cong \Hom_{\oK}(L_1, L_2)\otimes_{\oK} K.$$ 
Hence, 
$$\mathcal L_G(\Pi_1,\Pi_2)\cong \Hom_G( L_1, L_2)\otimes_{\oK} K\cong  \Hom_P(L_1, L_2)\otimes_{\oK} K\cong 
\mathcal L_P(\Pi_1,\Pi_2).$$
\end{proof} 

\begin{prop}\label{sub} Let $\pi$ be a smooth $\oK[G]$ module of finite length, and suppose that the irreducible subquotients of $\pi$ are either supersingular
or characters then every $P$-invariant $\oK$-submodule of $\pi$ is also $G$-invariant.
\end{prop}
\begin{proof} Let $\pi'$ be $\oK[P]$ submodule of $\pi$. If $\sigma$ is an irreducible subquotient of $\pi$ then by Theorem \ref{irred} $\sigma|_P$
is also irreducible, hence $\pi$ and $\pi'$ are $\oK[P]$ submodules of finite length.  

Let $\tau$ be an irreducible  $\oK[P]$-submodule of $\pi'$. Since $\pi$ is a finite length $\oK[G]$ module, the submodule 
$\langle G \centerdot \tau\rangle$ is of finite length. Let $\sigma$ be a $G$-irreducible quotient of $\langle G\centerdot \tau\rangle$.  
Since $\tau$ generates $\langle G\centerdot \tau \rangle$ as a $G$-representation, the $P$-equivariant composition:
$$ \tau\rightarrow\langle G\centerdot \tau\rangle\rightarrow \sigma$$ 
is non-zero, and since $\tau$ is irreducible, it is an injection. Now $\sigma|_P$ is irreducible, so the above composition 
is an isomorphism. Theorem \ref{restP} and Lemma \ref{char} imply that
$\tau$ is $G$-invariant and isomorphic to $\sigma$. By induction on the length of $\pi'$ as an $\oK[P]$-module
$\pi'/\tau$ is a $G$-invariant $\oK$-submodule of $\pi/\tau$. Since $\pi'$  is the set of elements of $\pi$  whose image in $\pi/\tau$ lies in $\pi'/\tau$,
 $\pi'$ is $G$-invariant.   
\end{proof}
 
\begin{cor} Let $\Pi$ be a  unitary $K$-Banach space representation of $G$, let $\lVert\centerdot\lVert$ 
be a $G$-invariant norm defining the topology on $\Pi$. Set
$$L=\{v\in \Pi: \lVert v\lVert\le 1\}.$$
Suppose that $L\otimes_{\oK}\Fbar$ is a finite length $\oK[G]$ module, and the 
irreducible subquotients are either supersingular or characters, then 
every closed $P$-invariant subspace of $\Pi$ is also $G$-invariant.
\end{cor}
\begin{proof} Let $\Pi_1$ be a closed $P$-invariant subspace of $\Pi$. Set $M=\Pi_1\cap L$, 
then $M$ is an open $P$-invariant lattice in $\Pi_1$.
Proposition \ref{sub} implies that for all $k\ge 1$, $M/\varpi_K^kM$ is a $G$-invariant 
$\oK$-submodule of $L/\varpi_K^k L$. By passing to the limit 
we obtain that $M$ is a $G$-invariant $\oK$-submodule of $L$. Since $\Pi_1=M\otimes_{\oK}K$ we obtain 
the claim. 
\end{proof}


\begin{thebibliography}{14}
\bibitem{bl}\textsc{L. Barthel and R. Livn\'e}, `Irreducible modular representations of $\GL_2$ of a local field', \textit{Duke Math. J.} 75 (1994)
261-292.
\bibitem{berger} \textsc{L. Berger} `Repres\'entations modulaires de $\GL_2(\Qp)$ et repr\'esentations galoisiennes de dimension 2', Preprint October 2005.
\bibitem{breuil}\textsc{C. Breuil},  `Sur quelques repres\'entations modulaires et $p$-adiques de $\GL_2(\QQ_p)$ I', \textit{Compositio Mathematica}
138 (2) (2003), 165-188.
\bibitem{p} \textsc{V. Paskunas}, `Coefficient systems and supersingular representations of $\GL_2(F)$', \textit{M\'emoires de la SMF} 99 (2004). 
\bibitem{schn} \textsc{P. Schneider}, \textit{Nonarchimedean Functional Analysis}, Springer-Verlag, 2001.
\bibitem{vig} \textsc{M.-F. Vign\'eras} `S\'erie principale modulo $p$   de groupes r\'eductifs $p$-adiques', Preprint April 2006.

\end{thebibliography}
\end{document}